\documentclass[12pt]{amsart}
\usepackage{pifont,geometry,latexsym,graphicx,mathptmx}
\usepackage{amssymb,amsmath,amsfonts,amscd,stmaryrd}

\newtheorem{Theorem}{Theorem} 

\newtheorem{Remark}[Theorem]{Remark}

\newtheorem{Definition}[Theorem]{Definition}
\newtheorem{Conjecture}[Theorem]{Conjecture}

\numberwithin{Theorem}{section}
\numberwithin{equation}{section}

\newcommand{\h}{\hspace*{.24in}}
\newcommand{\hhh}{\hspace*{.72in}}

\begin{document}
\title{On the conjecture about Morrey quasiconvexity in $L^{\infty}$}
\maketitle
\centerline{Hung Vinh Tran}
\centerline{\small Mathematics department - UC Berkeley}
\centerline{\small Email: tvhung@math.berkeley.edu}

\date{} 

\begin{abstract}
We study the difference between weak Morrey quasiconvexity and strong Morrey quasiconvexity in $L^{\infty}$. We point out some relations as well as give one example to show that weak Morrey quasiconvexity cannot imply strong Morrey quasiconvexity.
\end{abstract}
\section{Introduction}
In this short note, we try to understand about Morrey quasiconvexity in some $L^{\infty}$ variational problems. It is still an important question to characterize all sort of notions about quasiconvexity and to point out useful examples in applications.\\ \\
In \cite{BJW}, Barron, Jensen and Wang introduce the notion of strong Morrey quasiconvexity, which is the necessary and sufficient condition for the existence of solution of $L^{\infty}$ variational problem. They also introduce the notion of weak Morrey quasiconvexity, which is the analog of the usual Morrey quasiconvexity in $W^{1,p}$. It is known that in the scalar case (more precise statement will be introduced in the below section), those two definitions are equivalent. However, for the vector case, they don't know whether those two definitions are equivalent or not.\\\\
I show that in the vector case, those two definitions are not equivalent as in the below section. In fact, the example I introduce is quite simple and may not have any real application, but it indeed gives what we need. I would like to thank Hoang Tran for reading this note and giving useful suggestions. More on quasiconvexity and PDE can be found in the very nice book of Evans \cite{E1}.

\section{Main results}
Throughout this note, $Q$ will be denoted by the unit cube in the space $\mathbb R^{n}$ .\\
In \cite{BJW}, Barron, Jensen and Wang give the definitions of Morrey quasiconvexity as following
\begin{Definition}
A measurable function $f: \mathbb R^{nm} \to \mathbb R$ is said to be (strong) Morrey quasiconvex, if for any $\epsilon>0$, for any $A \in \mathbb R^{nm}$, and any $K>0$, there exists a $\delta = \delta(\epsilon, K,A)>0$ such that if $\varphi \in W^{1,\infty}(Q,\mathbb R^m)$ satisfies
\begin{equation}
||D \varphi||_{L^{\infty}} \le K,~ \max_{x \in \partial Q} |\varphi(x)| \le \delta,
\notag
\end{equation}
then
\begin{equation}
f(A) \le \mbox{ess~sup}_{x \in Q} f(A + D \varphi(x)) + \epsilon.
\notag
\end{equation}
\label{SMQ}
\end{Definition}
\begin{Definition}
A measurable function $f: \mathbb R^{nm} \to \mathbb R$ is said to be weak Morrey quasiconvex, or $(0,0)$ Morrey quasiconvex, if for any $A \in \mathbb R^{nm}$, and $\varphi \in W^{1,\infty}_{0}(Q, \mathbb R^m)$ we have
\begin{equation}
f(A) \le \mbox{ess~sup}_{x \in Q} f(A + D \varphi(x)).
\notag
\end{equation}
\label{WMQ}
\end{Definition}
It is obvious that strong Morrey quasiconvexity implies weak Morrey quasiconvexity.
\begin{Definition}
A measurable function $f:\mathbb R^{nm} \to \mathbb R$ is said to be quasiconvex if for any $s \in \mathbb R$, we have $f^{-1}((-\infty, s])$ is convex.
\label{Q}
\end{Definition}
For the case either $n=1$ or $m=1$, Barron, Jensen and Wang prove the following theorem
\begin{Theorem}
If either $n=1$ or $m=1$ and $f$ is lower semicontinuous, then the following notions are equivalent:\\\h
(i) $f$ is quasiconvex.\\\h
(ii) $f$ is polyquasiconvex.\\\h
(iii) $f$ is (strong) Morrey quasiconvex.\\\h
(iv) $f$ is weak Morrey quasiconvex.
\label{T1}
\end{Theorem}
\begin{Remark}
Here is one example of quasiconvex function when either $n=1$ or $m=1$.\\
In $\mathbb R^2$, take two arbitrary point $A$ and $B$ and denote by $AB$ the line segment between them.\\
 We define the function $f$ such that $f(x)=0$ for $x \in AB$ and $f(x)=1$ otherwise.\\
 It's easy to see that $f$ is lower semicontinuous and $f$ is quasiconvex. Hence $f$ is also strong, weak Morrey quasiconvex.
\label{E1}
\end{Remark}
It's still an open question that whether Morrey quasiconvexity and weak quasiconvexity are equivalent for $m,n>1$. In , Barron \cite{BJW}, Jensen and Wang give the following conjecture:
\begin{Conjecture}
For the case $m,n>1$ and $f$ is lower semicontinuous, strong Morrey quasiconvexity and weak Morrey quasiconvexity are not equivalent.
\label{C1}
\end{Conjecture}
Now, we prove this conjecture by point out one example of the function that is weak but not strong Morrey quasiconvex.\\
We consider for the case $m=n=2$.\\
Let $e_1,~e_2,~e_3,~e_4$ be the basis of $\mathbb R^4$ as normal. We denote by $M,~N,~R$ the three points in $\mathbb R^4$ as following $M=(1,0,0,0)$, $N=(0,0,1,0)$ and $R=(1,0,1,0)$ and $OM$, $ON$, $MR$, $NR$ be the four line segments.\\
Let $S$ be the square $OMRN$ ($S=OM \cup ON \cup MR \cup NR$).\\
Let $f: \mathbb R^4 \to \mathbb R$ such that $f(x)=0$ for $x \in S$ and $f(x)=1$ otherwise.\\
Firstly, it's easy to see that $f$ is lower semicontinuous and $f$ is not quasiconvex since $\{ x | f(x) \le 1/2\}$ is not convex.\\\\
We will prove that $f$ is weak Morrey quasiconvex.\\
To make things clear, we have several rules of writing as followings:\\
For $P \in \mathbb R^4$, we write $P=(P_1,P_2)$ for $P_1,~P_2 \in \mathbb R^2$.\\
For $\varphi \in W^{1,\infty}_{0}(Q,\mathbb R^2)$ we write $\varphi=(\varphi_1,\varphi_2)$ for $\varphi_1,~\varphi_2 \in W^{1,\infty}_{0}(Q,\mathbb R)$.\\
We need to prove:
\begin{equation}
f(P) \le \mbox{ess~sup}_{x \in Q} f(P + D \varphi(x)).
\notag
\end{equation}
or by writing in the new way:
\begin{equation}
f(P_1,P_2) \le \mbox{ess~sup}_{x \in Q} f(P_1 + D \varphi_1(x),P_2 + D \varphi_2(x)).
\label{Equ1}
\end{equation}
Let's consider two cases:\\\h
{\bf Case 1} If either $P_1 \notin OM$ or $P_2 \notin ON$. WLOG, we may assume that $P_1 \notin OM$.\\
Consider $g: \mathbb R^2 \to \mathbb R$ such that $g(x)=0$ for $x \in OM$ and $g(x)=1$ otherwise. As we already note in Remark $\ref{E1}$, g is weak Morrey quasiconvex. Hence,
\begin{equation}
1=g(P_1) \le \mbox{ess~sup}_{x \in Q} g(P_1 + D \varphi_1(x)).
\label{Equ2}
\end{equation}
So, we get the set $\{x \in Q| P_1 + D \varphi_1(x) \notin OM\}$ has positive Lebesgue measure or $|\{x \in Q| P_1 + D \varphi_1(x) \notin OM\}| >0$. Furthermore, we have:
\begin{equation}
\{x \in Q| P_1 + D \varphi_1(x) \notin OM\} \subset \{x \in Q| f(P+ D \varphi (x)) = 1\}.
\label{Equ3}
\end{equation}
Therefore, we get:
\begin{equation}
|\{x \in Q| f(P+ D \varphi (x)) = 1\}|>0.
\label{Equ4}
\end{equation}
So,
\begin{equation}
\mbox{ess~sup}_{x \in Q} f(P + D \varphi(x))=1 \ge f(P),
\notag
\end{equation}
we get the result.\\\h
{\bf Case 2} If both $P_1 \in OM$ and $P_2 \in ON$. If $\varphi=0$ then we obviously get the result. Hence, assume that $\varphi \ne 0$. WLOG, we may assume $\varphi_1 \ne 0$.
We will prove the similar result as above:
\begin{equation}
|\{x \in Q| P_1 + D \varphi_1(x) \notin OM\}| >0,
\notag
\end{equation}
and then get the result.\\
Suppose not, then $P_1 + D \varphi_1(x) \in OM$ a.e. This implies $\dfrac {\partial \varphi_1}{\partial x_2}=0$ a.e.\\
Since $\varphi_1|_{\partial Q} = 0$, we then get $\varphi_1=0$ a.e., which is absurd.\\\\
Lastly, we will prove that $f$ is not strong Morrey quasiconvex.\\
Take the point $P=(1/2,0,1/2,0) \in \mathbb R^4$. Notice that $f(P)=1$.\\
Let $\phi:\mathbb R \to \mathbb R$ be the periodic "zig-zag" function defined by
\begin{equation}
\phi(t)=
\left\{ \begin{aligned}
t \hhh & \mbox{if}~0 \le t \le 1/2 \vspace{.05in}\\
1-t \hhh & \mbox{if}~1/2 \le t \le 1,\\
\end{aligned} \right. 
\notag
\end{equation}
and $\phi(t+1)=\phi(t)$ for $t \in \mathbb R$.\\
For any $\epsilon>0$,  we then define $\varphi_1(x)=\varphi_2(x) = \dfrac{\epsilon}{2} \phi(\dfrac{x_1}{\epsilon})$.\\
Note that we have:  $\dfrac{\partial \varphi_1}{\partial x_2} = \dfrac{\partial \varphi_2}{\partial x_2}=0$ and $|\dfrac{\partial \varphi_1}{\partial x_1}| =|\dfrac{\partial \varphi_2}{\partial x_1}| =\dfrac{1}{2}$ a.e.\\
Hence $ P + D \varphi \in S$ a.e. Therefore, $\mbox{ess~sup}_{x \in Q} f(P + D \varphi(x))=0$.\\
We get $f$ is not strong Morrey quasiconvex.

\begin{Remark}
This example is sort of trivial. The function $f$ is equal to $1$ a.e.  However, it's indeed the example to show that the two notions are not equivalent. I don't know if $f$ makes any sense in applications. I think it is quite important to warn that in this problem, $f$ is required to be uniquely defined at each point. Or more clearly, $f$ cannot be defined just almost everywhere or $f$ cannot be defined in such class of $L^{\infty}$ or other similar classes.
\label{R1}
\end{Remark}


\begin{thebibliography}{xx}
 \bibitem{BJW} {E.N. Barron, R.R. Jensen, C.Y. Wang},
{\it    Lower semicontinuity of $L^{\infty}$ functionals}, Ann. I. Pointcare-AN 18,  {\bf 4} (2001), 495-517.
 \bibitem{E1} {L. C. Evans},
{\it    Weak convergence methods for nonlinear partial differential equations}, CBMS 74, American Mathematical Society, 1990.
\end{thebibliography}
\end{document}